\documentclass[11pt]{amsart}
\usepackage[T1]{fontenc}
\usepackage[ansinew]{inputenc}
\usepackage{amssymb,amsmath,amsthm,eucal,hyperref,mathrsfs,array,graphicx,xcolor,fullpage}
\usepackage[all]{xy}

\theoremstyle{definition}

\def \eps {\epsilon} 

\def \N {\mathbb{N}}

\def \R {\mathbb{R}} 

\def \Z {\mathbb{Z}}

\title{A simple renormalization flow for FK-percolation models}
\author{Wendelin Werner}
\address{D-Math, ETH Z\"urich, R\"amistr. 101, 8092 Z\"urich, Switzerland.}
\email{wendelin.werner@math.ethz.ch}

\begin{document}

\begin {abstract}
We present a  setup that enables to define in a concrete way a renormalization flow for the FK-percolation models from statistical physics (that are 
closely related to Ising and Potts models). In this setting that is applicable in any dimension of space, one can interpret perturbations of the critical (conjectural) scaling limits
in terms of stationary distributions for rather simple Markov processes on spaces of abstract discrete weighted graphs. 
\end {abstract}

\maketitle


\section {Background} 
\noindent
{\bf General introduction.}
The study of critical phenomena and the fine understanding of phase transitions have been one of the success stories of theoretical physics in the last seventy years. Many ideas and approaches to these questions have been proposed, from explicit combinatorial identities to renormalization group or field theory, a number of which have in turn given rise to  important developments in several branches of mathematics.

In particular, the renormalization group has turned out to be a powerful, versatile and remarkably successful idea that enabled theoretical physicists to derive and predict numerous features related to phase transition and critical phenomena for physical systems of various types and in several different dimensions of space. The bibliography on this topic is gigantic (one can mention the names of Fischer, Wilson, Kadanoff, Symanzik, De Gennes and many others, see for instance the reference lists in \cite {Ca,Z}). 
A rough basic idea in the context of models from statistical physics is that when considered at their critical point, their behavior can
 remain random at all scales, and becomes in fact statistically self-similar in the scaling limit. This suggests that the continuous limit can be viewed as a fixed point of a certain renormalization map (that goes from one scale to the next one), or of a renormalization flow (when one zooms in or out continuously). 
 
 If for a given map (corresponding to a given model, say for instance the Ising model) and a given dimension of space, there would exist a unique fixed point for this renormalization map, then it would explain why this model when taken on different lattices (and at their critical point) behaves always in the same way in the scaling limit. This phenomenon that is often referred to as ``universality'' is today still far from being mathematically well-understood, even if a lot of progress has been made on some specific models in two-dimensions. 
 
 The simplest example of such a renormalization map can easily be heuristically described in the case of planar percolation. 
 Suppose for instance that one has a way to color at random (i.e., according to some probability measure $P$) the unit square $[0,1]^2$ into two colors (white and black) and
 considers the collection of all all-white paths.
Then, one samples four independent copies of this coloring, one puts them next to each other in order to form a larger square of side-length two, and scales down the picture by a factor of $2$. In this way, we have now obtained a new random configuration in the unit square that is distributed according to some new distribution $Q(P)$ (in this new configuration, clusters from different squares can be concatenated in order to form a larger one). 

There are two obvious trivial fixed points for the map $Q$: The probability measures that colors the square all black almost surely, or all white almost surely. 
The question is whether there are any other fixed points of this renormalization map $Q$ ie. non-trivial distributions $P$ such that $Q(P)=P$. 
The universality  conjecture for percolation is then that there exists a unique non-trivial fixed point for $Q$, which would describe the scaling limit of all planar critical percolation models. 
There are technical difficulties in making the state space in which $P$ is defined precise (here one for instance better views the coloring as a collection of all-white paths), and also how to properly define $Q$, but many of these has been overcome  by Schramm and Smirnov in their paper \cite {SchSm}, where they have been able to exhibit one such a non-trivial fixed point (this fixed point is also related to Tsirelson's theory of black noise, \cite {T}). The uniqueness question remains open. 

There is however a conceptual problem when one tries to make sense of such a renormalization map for dependent models, such as the Ising model. Indeed, one cannot just glue together configurations chosen independently according to $P$ in different squares, as they are not independent anymore; in fact,  they depend on each other in a very complex way in the large-scale limit, and understanding this dependence is precisely the heart of the matter in the study of those models. It appears therefore very difficult to generalize the previous percolation setting to other models.

\medbreak
\noindent
{\bf The FK percolation model.}
The FK-percolation model (sometimes also called the random cluster model) named after Fortuin and Kasteleyn \cite {FK} who introduced it around 1970 is a classical model of statistical mechanics that is very closely related to the Ising and Potts models (see also \cite {GG} and the references therein). 

When $G$ is a finite graph (we will denote by $V$ its set of vertices and by $E$ its set of unoriented edges), FK-percolations form a  family of distributions on the set of functions $w$ from $E$ into $ \{  0, 1 \}$. When $w(e)=0$, then the edge $e$ is said to be closed, while when $w (e)=1$, it is declared open. Hence, such a map $w$ defines in fact a subgraph $G(w):=(V,E_w)$ of $G$, where $E_w \subset E$ is the set of open edges for $w$.
We denote for each $w$ by $o(w)$ the number of open edges for $w$, and by $k(w)$ the number of connected components of $G(w)$.  
Of course, one has $k(w) \le \# V$ and the number of closed edges for $w$ is $\# E - o(w)$. 

For each $q >0$ and each $p \in [0,1]$, we define $\pi = \pi (p) = p / (1-p)$ and the FK-percolation distribution $P_{p,q}=P_{p,q,G}$ is defined by 
$$ 
P_{p,q} (w) =  \frac 1 {Z_{p,q,G}}  \times \pi^{o (w)} q^{k(w)},$$
where $Z_{p,q,G}$ is just the renormalizing constant (usually called the partition function) defined by
$$ 
Z_{p,q,G} = \sum_w \pi^{o(w)} q^{k(w)}
$$
so that $P_{p,q}$ is indeed a probability measure. 

In the case where $q=1$, $P_{p,q,G}$ is just the Bernoulli percolation measure, where each edge $e$ is independently declared open or closed with probability $p$ and $1-p$. 
In the other cases, one gives an additional weight proportional to $q^{k(w)}$ to each configuration. When $q$ is a positive integer, this FK-measure is 
very directly related to the $q$-state Potts model (which is equal the Ising model when $q=2$): The connectivity properties of $w$ are then describing the correlations for the Potts/Ising model, see again for instance \cite {GG}, which is one of the main motivations to study this FK-percolation model and its phase transition.

When $q \not=1$ and $G$ is a connected graph, then the state of the edges are not independent anymore. However, it is immediate to check that for each $q \ge 1$, and each edge, one always has $ p/q \le P_{p,q} (w(e) =1 ) \le p$.

For a given graph and given $q$, it is also possible to choose a different value $p_e$ for each edge $e$. For every function $\bar p = (p_e, e \in E)$ from $E$ into $[0,1]$, one  defines the corresponding FK measure via 
$$ P_{\bar p,q} (w) = \frac 1 {Z_{\bar p,q,G}}  \left( \prod_{e \in E} (\pi_e 1_{w(e)=1} + 1_{w(e)=0})  \right) q^{k(w)}$$
where $\pi_e = \pi (p_e) = p_e / (1-p_e)$ and $Z_{\bar p, q, G}$ again denotes the renormalization constant so that this is a probability measure.

An important and useful feature when $q \ge 1$ is that increasing events are positively correlated (this is the so-called FKG inequality) meaning in for instance that for two different edges $e$ and $e'$, one has 
$$ P_{\bar p,q,G} ( w(e) = w (e') = 1 ) \ge P_{\bar p,q,G} ( w(e)=1 ) P_{\bar p,q, G} (w(e') = 1). $$
More generally, the FKG inequality enables to couple (for a fixed $p$ and $q$) the FK-percolation measures on two graphs $G$ and $G'$ when $G \subset G'$ in such a way that almost surely,  any edge of $G$ that is open in the FK-realisation in $G$ remains open for the FK-realisation in $G'$. 
This shows that if one considers 
an increasing sequence $G_n$ of finite graphs, one can define the probability measure $P_{p,q,G}$ where $G$ is the limit (i.e. union) of the graphs $G_n$ as the 
limit of the measures $P_{p,q,G_n}$ (see for instance \cite {GG}) -- note that this will work as soon as the vertex set of $G$ is countable and that a vertex is allowed to have 
infinitely neighbors in the graph, which will be the case in the present paper. This way to make sense of FK-measures on infinite graphs $G$ works also for non-constant functions $\bar p$ (defined on the edge-set of $G$).

Let us briefly summarize some well-known features of the FK-percolation model and related conjectures regarding its phase transition. The idea is to fix $q \ge 1$, and to study the measures $P_{p,q}$ on an infinite $d$-dimensional transitive graph (like $\Z^d$) and their features when one lets $p$ vary.  
\begin {itemize} 
 \item In a graph like $\Z^d$ for $d \ge 2$, the measures $P_{p,q}$ exhibits a phase transition: There exists a critical value $p_c = p_c (q)$ (depending on the graph) in $(0,1)$ such that for $p < p_c$, a configuration chosen according to 
 $P_{p,q}$ has no infinite cluster, while for $p > p_c$, it almost surely has a (unique) infinite cluster. In the case $q =2$, this exactly corresponds to the phase transition of the Ising model (that models the ferromagnetic phase transition of iron at the Curie point).   
 \item For a given $d$, when $q$ is not too large, the critical measure $P_{p_c, q}$ is supposed to exhibit interesting scale-invariance features, and it is conjectured that one can define a continuous scaling limit in some appropriate configuration space (when the dimension of the space is large, above the so-called critical dimension for the given value of $q$, the model at criticality does not exhibit such scale-invariance properties and the phase transition is then said to be discontinuous because the map $p \mapsto 
 P_{p,q} (0 \hbox { is in an infinite component})$ is discontinuous at $p_c$). For those scale-invariant models, one expects polynomial decay of certain quantities governed by so-called critical exponents. For instance, 
 the probability that one can find an open path from the origin to the boundary of a disc of size $R$ around the origin should decay like $R^{-\alpha + o(1)}$ as $R \to \infty$. Similarly, the probability that two given neighboring sites belong to two different connected components of $w$ that each have a diameter at least $R$ should decay like 
 $R^{-\beta + o(1)}$ as $R \to \infty$, where $\alpha$ and $\beta$ are (dimension-dependent but lattice-independent) constants. The exponent $\beta$ is important as it describes the probability for an edge to be pivotal (if one just changes its state from closed to open, one creates a much larger cluster). 
 \item The {\em universality  conjecture} is that for a given $q$ and a given dimension of space (that is below a critical dimension associated to $q$), the critical models defined on different $d$-dimensional grids do all converge to the same scaling limit. As mentioned above, the rule of thumb is that this scaling limit is the unique fixed-point (in $d$ dimensions) of a certain renormalization map, that should reflect the asymptotic scale-invariance and the fact that the lattice effects disappear. The fact that for a given $d$ and a given $q$, the exponent $\alpha$ is the same for all $d$-dimensional grids would be a side-effect of this universality. 
\end {itemize}

\medbreak 
\noindent
{\bf Near-critical percolation, the work of Garban, Pete and Schramm.}
The phase transition for two-dimensional percolation in the case of site percolation on the triangular lattice is now very well-understood. The large-scale behavior of the critical 
model and its continuous scaling limit has been derived thanks to Smirnov's discrete holomorphicity approach \cite {Sm}, and its scaling limit can be described via Schramm's SLE$_6$ \cite {S} process, which in turn allows to determine the values of the exponents such as $\alpha$ and $\beta$ that turn out to be equal to $5/48$ and $5/4$ respectively (see \cite {LSW5,SW}). Furthermore, as explained above,  it can be viewed as the fixed point of a renormalization map (in some appropriately defined probability space) \cite {SchSm}. A related question, that will be relevant to the approach that we will develop in the present paper in the case of FK-percolation, is to describe the behavior of percolation configurations when the probability for a site to be open is very close, but not equal, to the critical value.

To treat this question in the case of usual percolation, it is natural to couple the percolation configurations for all values of $p$ on one single probability space, by sampling for each site $x$ of the triangular lattice ${\mathcal T}$ an independent uniformly distributed random variable $X(x)$ on the interval $[0,1]$, and to define then for each $p$, the configuration $u_p = (u_p (x), x \in {\mathcal T})$ with $u_p (x) = 1_{X(x) \le p }$. This provides a natural coupling of percolations with all parameters $p$, and if one looks at the map $p \mapsto u_p$, one basically observes the movie where all sites become progrssively open in a random uniformly chosen order. 

When $p = p_c - \epsilon$ for a small $\epsilon$, the configuration $u_p$ can be interpreted as the configuration $u_{p_c}$ where a small uniformly chosen fraction of the edges that were open for $u_{p_c}$ have been closed. The scenario that has been proved to hold in the sequence of papers by Garban, Pete and Schramm \cite {GPS1,GPS2,GPS3}, is the following: 

\begin {itemize} 
 \item When one chooses $N = N (\eps)$ appropriately (or conversely, $\eps = \eps (N)$) and rescales the picture by $1/N$, the obtained random picture (in fact rigorously defined  in terms of its macroscopic connections) converges to a non-trivial limit (when $\eps \to 0$ and $N \to \infty$). This limit is not scale-invariant: When one zooms in on  a smaller scale and then scales it up, the obtained picture corresponds to 
 changing the factor $\epsilon$ by some constant factor i.e. to open some edges uniformly at random. 
 \item This limit can be described by first sampling the scaling limit of critical percolation, and then on this continuous picture, removing a Poissonian collection 
 of pivotal points (where one cuts a macroscopic connection when one removes this point) sampled according to an intensity measure on pivotal points, that can be read off from the continuous critical picture.  
\end {itemize}
In other words, one can pass the ``sample the critical model and erase some of its sites/edges at random'' procedure to the continuous 
limit. This short summary covers in fact a number of highly non-trivial facts and we refer to \cite {GPS1,GPS2,GPS3} for background and details. 

This shows in particular that this near-critical scaling limit is in fact invariant under a renormalization map that can be loosely described as follows: 
\begin {itemize} 
 \item Consider four independent copies of the near-critical percolation picture in four given squares. Put them next to each other, look at the obtained picture, and then scale it down by a factor $2$ (this first step is exactly like for the critical percolation renormalization). 
 \item Then, open uniformly at random (i.e. in a Poissonian way, according to a well-chosen intensity measure) some of the points between adjacent clusters (and thereby creating new connections). 
\end {itemize}
This is exactly the type of Markovian dynamics (the random opening of some edges/points corresponds to a random dynamics) that we will use for other FK-models.

\medbreak 
\noindent
{\bf Renormalization and near-critical models.}
 The ``fixed points'' that appear via renormalization group analysis are of two types: The ones that describe the scaling limit of the system exactly at the critical point, and the ones that describe the ``near-critical'' behaviour i.e., how the system behaves in the limit where the size of the system and the parameter (typically the temperature, so we will denote it by $T$ in this paragraph) are tuned (via what is often referred to as finite-size scaling) so that in the large scale limit, the system is neither the same as at its critical point, nor trivial. Typically, just as in the near-critical percolation case, one fixes a parameter $\lambda$, and considers the limit when $N \to \infty$ where $N$ denotes the size of the system) at temperature $T = T_c + \lambda N^{-\beta}$ for some well-chosen critical exponent $\beta$. This near-critical scaling limit is then not scale-invariant anymore, but a change of scale just corresponds to a change of the parameter $\lambda$. These near-critical scaling limits are 
intrinsically related to the critical ones (they for instance provide the same critical exponents etc.) even though they have also essential differences (see for instance \cite {NW}). 
The framework that we will develop in the present paper will be more related with this second (i.e. near-critical) approach. 

 Note however that there are more than one parameter that one can play with in order to perturb the discrete model. For instance, in the case of the Ising model, there is the inverse temperature or the external magnetic field (and in fact many other possible perturbations). So, the critical model appears to be critical in ``many directions'' i.e. it is a renormalization fixed point in a parameter space with many dimensions. So, to each critical model, one can associate many different near-critical perturbations (continuous scaling limts) using the finite-size scaling ideas mentioned above; there is basically one for each ``perturbation direction'' in the parameter space. 
 
 As it turns out, the dynamics of the phase transition with respect to the (inverse) temperature is quite intricate for FK-percolation when $q \not= 1$. While it is possible (via FKG-inequality-type arguments) to couple in a fairly natural way the different FK models for different values of $p$ in an increasing manner, the way in which edges appear one after the other in this coupling is highly correlated and still not well-understood, as pointed out in \cite {DGP}.
 
 The originality of the approach developped in the present paper is maybe its choice of direction that enables to encapsulate the renormalization operation in a fairly simple way. The basic idea is the following: We let the edges of a given FK-percolation model (i.e. for a given choice of $p$ and $q$) appear one by one in continuous time and in uniformly chosen order. This process turns out to have very nice Markovian properties that we shall now describe and exploit.

\section {A setup for renormalization}
\noindent
{\bf A simple observation.} 
Suppose that we are given a finite graph $G$, and a sample $w$ of the FK-percolation measure $P_{p,q,G}$ for $q \ge 1$. We also choose a value $\alpha \in (0,1)$. Conditionally on $w$, we sample another configuration $v$ by simply tossing an independent $\alpha$-coin for each edge that is open for $w$ in order to decide whether this edge is still open for $v$ (with probability $\alpha$) or closed (with probability $1-\alpha$); all edges that are closed for $w$ remain closed for $v$, so that $v \le w$.  The joint distribution of $(v,w)$ is then 
$$ 
P ( v, w) = \frac {1}{Z_{p,q,G}} \times  \pi^{o(w)} q^{k(w)} \alpha^{o(v)}  (1-\alpha)^{o(w) - o(v)} 
 = \left[ \frac {q^{k(v)} (\alpha \pi)^{o (v)} }{Z_{p,q,G}} \right] \times (\pi (1-\alpha))^{o(w) - o(v)} q^{k(w)- k(v)} 
.$$ 
It follows that if we define the new graph $G/v$ obtained from $G$ by collapsing all edges that are open for $v$ (i.e., we identify the two end-points $x$ and $y$ of all the edges that are open for $v$, and then erase this edge as well as all other edges that join $x$ and $y$ -- this graph should not be confused with the graph $G (v)$ defined earlier), then one can view $w-v$ as a configuration on this graph, and the previous expression shows that {\em conditionally on $v$, $w-v$ is distributed exactly like an FK-percolation configuration on $G/v$ with parameters $p'$ and $q$, where 
$\pi (p') = \pi (p)  \times (1-\alpha)$.
}

Note that when $\alpha = \epsilon$ is very close to $0$, at first order i.e. with a probability $1 - O(\epsilon)$, $v$ has no open  edge, and
$ \pi (p') = \pi(p) e^{- \epsilon} + o (\epsilon)$.  The probability that the edge $e$ is open for $v$ is $ \epsilon P_{p,q, G} ( w(e)=1)  + o (\epsilon)$, while the probability that at least two different edges are open for $v$ is of the order $\epsilon^2$. 

The very same observation works in the case of non-constant edge-weights $\bar p = (p_e, e \in E)$: One just needs to choose $\bar p'$ with $\pi (p_e') = \pi (p_e) \times (1-\alpha)$ in order to define the conditional distribution of $w$ given $v$. 

\medbreak

\noindent
{\bf A slightly different setup.}
This leads naturally to the following equivalent description of FK-percolation. We are now considering {\em weighted graphs} i.e. couples $W=(V,c)$, where $V$ is a finite set of vertices, and $c$ is a function from $V \times V$ into $\R$ ($c(x,y)$ can be viewed as a conductance of the edge between $x$ and $y$, hence the letter $c$) with $c(x,y)=c(y,x)$ and $c(x,x)=0$. 
We will sometimes associate to $W$ the (non-weighted) graph $G = (V,E)$ where $E$ is the set of edges with non-zero $c$-weight.

We can then define for each $q$, the FK-model on this weighted graph (without specifying the value of the function $\bar p$, because this value will be implicitely determined by the function $c$) as the $P_{\bar p, q}$ FK-measure on $G$ (or on the complete graph with vertex set $V$), where the value of $p_e$ associated to the edge $e$ will be given by the formula $p_e = 1- e^{-c (e)}$. In other words, as $\pi (p_e) = e^{c(e)} - 1$, the distribution of $w$ is then 
$$ P_{W,q} (w) = \frac {1}{Z_{W,q}} \times q^{k(w)} \times  \prod_e \left ( 1_{w(e) =0 } +  (e^{c(e)} - 1 ) 1_{w(e) = 1} \right) 
= \frac {1}{Z_{W,q}} \times q^{k(w)} \times  \prod_{e : w(e) =1}  (e^{c(e)} - 1 ) $$
(here by slight abuse of notation, we write $Z_{W,q} := Z_{\bar p, q, G}$).  
Note that when $c(e)=0$, then the edge $e$ is anyway closed (whereas in the limit $c(e) = \infty$, it is anyway open). We note also that the probability that a given edge is open is a non-decreasing function of the function $c( \cdot)$ (via the FKG inequality). 

We then define for each edge $e$ (i.e. each pair of points in $V$), an independent exponential random variable $\xi (e)$ with mean $1$, that is also independent from $w$, and we then define for each $t \ge 0$,  the configuration $w_t 
= (w_t (e), e \in E) \in \{0,1\}^E$ as 
$$ w_t (e) = w(e) 1_{\xi (e) \le t}.$$
Note that at time $0$, all edges are closed for $w_0$, while almost surely for $t$ large enough, one has $w_t = w$. 

In the case where we start with all $p_e$'s equal to each other, the joint distribution of $(w_t,w)$ (for fixed $t$) is exactly that of the couple $(v,w)$ described above, with $\alpha = e^{-t}$. To each $t$, we can then associate to $w_t$ the weighted graph $W_t= (V_t, c_t)$ that is obtained by collapsing each of the connected components created by $w_t$ into a single site as described above (this defines $V_t$), and choosing for each two neighboring connected components $C$ and $C'$ a weight $c_t(C,C')$ equal to the sum of all $c(e)$'s where $e$ spans the set of edges of $W$ that join $C$ and $C'$. We will denote by $e_t$ the edges of $W_t$ (ie. of pairs of points in $V_t$ with non-zero edge-weight $c_t$). 

Now, the previous observation yields readily that the process $t \mapsto W_t$ is Markovian, and that its dynamics can be described via the two following rules: 
\begin {itemize}
\item The edge-weights have a erosion so that each of the values $\exp (c_t(e_t)) - 1$ decrease exponentially in time. 
In other words, when one starts the Markov dynamics at time $0$ with the weighted graph $W=(V,c)$, then during some time the vertex set $V_t$ does not change, but the weights evolve according to the rule 
$$  e^{c_t ( \cdot) }  - 1 = e^{-t} ( e^{c_0 (\cdot)} - 1 ) \quad ie. \quad {c_t( \cdot ) }  = \log (  1 + e^{-t} ( e^{c_0 (\cdot)} - 1 )).$$ 
This happens until the first edge collapses (if this happens at all) according to the rule described below. 
\item At time $t$, an edge $e_t$ of $W_t$ can be opened at a rate equal to the $P_{W_t,q}$-probability that it is open. When this happens, this corresponds to a jump of the Markov process to a new state $W_t'=(V_t', c_t')$, where the two extremities $x$ and $y$ of $e_t$ are collapsed into a single site $xy$, and where $c_t' ( xy, z) := c_t (x,z) + c_t (y,z)$ for all $z \in V_t \setminus \{ x,y\}$. 
\end {itemize}
\begin{figure}[ht!]
\centering
\includegraphics[width=10cm]{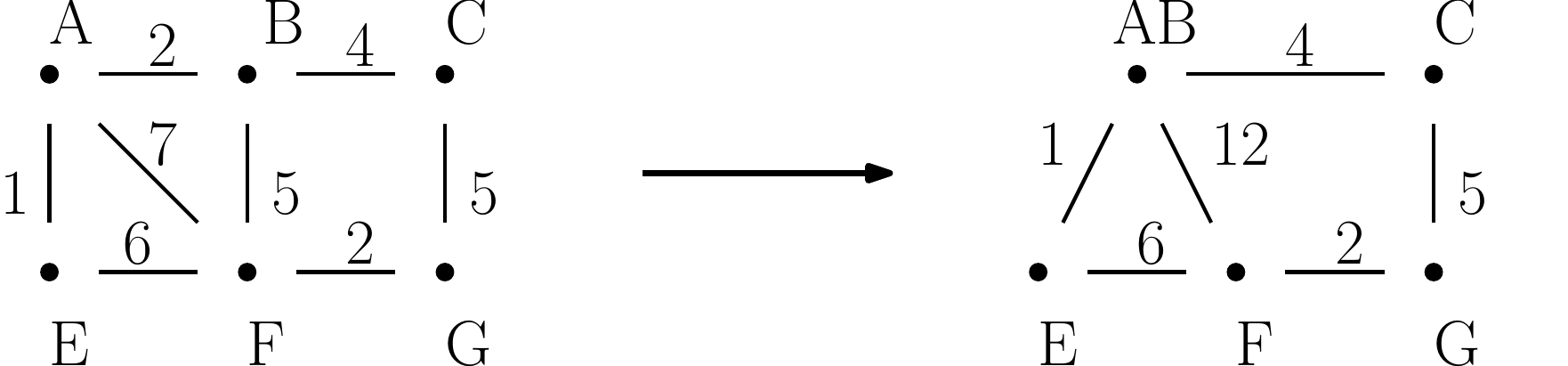} 
\caption{A jump of the process: Collapsing of the edge between A and B.}
\end{figure}

This is a fairly simple Markov process (that depends solely on the value of $q$) on the state of all finite weighted graphs. Note that it is clear from our definition and the interpretation in terms of $w_t (e) = w(e) 1_{\xi (e) \le t }$, that almost surely as $t \to \infty$, all $c_t (e_t)$ will eventually decay exponentially fast to $0$, and that the graph $W_t$ will stabilize to the graph obtained by collapsing all edges of $w$. 

\medbreak
\noindent
{\bf The infinite volume version of the Markov processes.}
It is not very difficult to extend the definition of each of these Markov processes (recall that there is one process for each value of $q$) to the state-space of infinite countable graphs, building on the existence of the infinite volume FK measure: Suppose 
that a weighted graph $W=(V,c)$ is given, where $V$ is countable and $c$ is some map from $V \times V$ into $\R_+$ (with the same conditions $c(x,y)=c(y,x)$  and $c(x,x)=0$ as above).
Sample the (infinite volume, with free boundary) weighted FK model on this infinite graph (note that the probability that a given edge $e=(x,y)$ is open for this configuration is the limit when $n \to \infty$ of the probability that it is 
open for the FK model on a finite approximation $W_n = (V_n, c)$ of the infinite graph if $V_n$ is increasing with $\cup V_n = V$). 
Then, one can define the dynamics using an auxiliary exponential random variable $\xi (e)$ for each edge as before. The obtained local dynamics (near a given site or edge) is then easily checked to be the $n \to \infty$ limit of the local dynamics when applied to the graphs $W_n$ (this is again due to the FKG inequality).  

Whether the initial graph (i.e. the initial values chosen for $c$) is subcritical or supercritical for the FK-model can then be read off from the behavior of the Markov process (in the former case, the weight of all edges will eventually go exponentially fast to $0$, and locally, the graph freezes to some final finite configuration, while in the latter case, a special vertex (i.e. a cluster) with infinitely many neighbors will appear somewhere (at finite distance of any given site) after some finite time, and this cluster will keep swallowing other points throughout the dynamics (after any given large time, it will still merge with some other clusters).

\medbreak
\noindent
{\bf Some rather trivial general comments on random infinite graphs.}
Suppose that two infinite weighted graphs $W=(V,c)$ and $W'=(V', c')$ are given. We say that they are equivalent (and write $W \sim W'$) if there exists a bijection $\varphi$ from $V$ onto $V'$ such that for all $x,y$ in $V$,
$c' (\varphi (x), \varphi (y))= c(x,y)$. Then, clearly, the previously described Markov dynamics applied to $V$ and applied to $V'$ preserve the equivalence i.e. one can couple the process started from $(V,c)$ with the process
started from $(V', c')$ in such a way that they remain equivalent at all times (i.e. $W_t \sim W_t'$). In other words, the Markov process does not depend on the actual labeling of the graph.

When one wants to {\em define} a random infinite graph (with countably many vertices, and where vertices can have countably many neighbors), one has to specify a space in which $V$ lives i.e. an embedding of $V$ into some given space. This embedding can be deterministic (i.e. $V$ can be a deterministic set, such as $\N$, which corresponds to an explicit enumeration of the vertex set) or random (the set $V$ can be a point process in some larger space such as $\R^d$ or the space of compact subsets of $\R^d$). 
If one is given any embedding, one can then define such an enumeration when the graph $W$ is connected by first choosing a root in the vertex set that will then be 
the first vertex, and then using some lexicographic-type rule in order to order the elements of $V$ (for instance explore them one by one starting from the root, using some rule involving the weights in order to decide on the ordering). A rerooting procedure would then correspond to applying a bijection $\varphi$ to the labels and it defines of course an equivalent graph.

When one is given a rooted weighted graph, and runs our Markovian dynamics, then one can view the graph $W_t$ as being rooted at the cluster than contains the initial root of $W_0$. 
For instance, when one starts with the graph $\Z^d$ with constant weights (or any other transitive connected weighted graph) and runs the Markovian dynamics, one can root it at time $0$ at some given point (say, the origin), and then root the graph $W_t$ at time $t$ on the cluster that contains the origin.
However, this enumeration procedures turn out to be not so convenient when one wishes to describe the asymptotic behavior of critical models (for instance, in the limit, the root ends up almost surely to be disconnected from all other sites as the cluster containing the origin freezes, so that the weights of all edges that touches this root go to $0$); this corresponds to the fact that one wants to look at the global behavior of the graph rather than at its local behavior near a given point, which is reminiscent of the difference between the local configuration of critical percolation on $\Z^d$ and its continuous scaling limit. 

A natural simple way to label the clusters at time $t$ is just to label them by the set of original vertices that they contain. This  carries of course a lot of information that is not needed in order to run the dynamics, and it becomes quite messy.  Since the  goal of this Markov process approach is precisely to get rid of this complexity when describing the fixed point that we have in mind, this is not so easy either. One can however simplifiy things a little by assign to each cluster a (different) point in $\R^d$ (we then say that 
the graph $W$ is embedded in $\R^d$).  For instance, when one starts the dynamics with $V = \Z^d$, when applying the Markov dynamics, one can apply a deterministic rule to decide how to label the cluster obtained from the merging of two clusters (for instance the label of the one of the two sites with largest boundary). This embedding approach provides one possible way to define a class of $d$-dimensional graphs:

We say that a probability measure defining a random weighted graph $G=(V,c)$ is a random translation-invariant $d$-dimensional graph if it satisfies:
\begin {itemize} 
 \item Translation-invariance: The vertex-set $V$ is a random subset of $\R^d$, and the obtained graph is translation-invariant (in distribution):  For all fixed $z \in \R^d$, if one defines $c_z ( \cdot, \cdot) = c ( \cdot - z, \cdot -z)$, then the laws of $(V,c)$ and $(V+z, c_z)$ are identical. 
 \item Local finiteness: For each $\eps >0$, the mean number of vertices $x \in V$ that are in the unit ball and have at least one neighbor with $c(x,y) > \eps$ is finite. 
 \item The graph is almost surely connected, and for each site $x$, the sum $\sum_y c(x,y)$ is finite (this quantity can be referred to as the perimeter of $x$). 
\end {itemize}

An example of such a graph is for instance $u+ \Z^d$, where $u$ is a uniformly chosen point in $[0,1]^d$. 
We can note that if one starts the dynamics from a transitive lattice such as $u+ \Z^d$ and constant (non-supercritical) $c$, then at any finite time $t$, the couple $(G_t, c_t)$ still is a random translation-invariant $d$-dimensional graph.
 
 Note that there exists alternative ways to try to define its dimension (for instance using isoperimetric properties of the graph $G^\eps$ obtained by keeping only the edges of weight greater than some small $\eps$) that we will not discuss here. Also, using ergodic properties of this $G^\eps$ can lead to a nice intrinsic definition of its rooting (for instance and loosely speaking: root it on the position of a random walk on $G^\eps$ after long time).

\medbreak
\noindent
{\bf The critical cases (heuristics and conjectures).}
Let us try to heuristically describe the conjectural asymptotic behavior of the Markov process when one starts with a critical value of the weights (and in a low enough dimension, so that this critical model exhibits asymptotic scale-invariance properties) in some finite-dimensional lattice.

Let us first focus on the very special case of percolation ie. $q=1$. Then, the state of the process after time $t$ corresponds to percolation with parameter $p_c (1-e^{-t})$, and it is therefore described in the two-dimensional case by the near-critical regime studied in \cite {GPS1,GPS2,GPS3}. In particular, modulo rescaling of space, the law of the merging of the (rescaled) clusters after time $t$ is conjecturally asymptotically independent of $t$. This suggests that the law of the pair $(G_t,c_t)$ does in fact converge to a stationary distribution, which can be obtained as in \cite {GPS2,GPS3} by first sampling the critical percolation scaling limit, then removing a Poisson point process of pivotal points (in the two-dimensional case, with respect to the four-arm measure $\mu$), and finally defining the collapsed graph (each site of the graph corresponds to a cluster, and the exponential of minus the edge-weight between two neighboring clusters is given by its pivotal measure mass). In this setting, we can for instance simply label the sites of 
the graph $G_t$ as the set of all the clusters, and when two sites collapse in the Markovian evolution, it corresponds to the reunion of two clusters into one (the label of the new site being now the larger cluster).
When one lets $t$ tend to $\infty$, one therefore ends up in the near-critical regime (without having to tune the near-critical window, the process loosely speaking self-adjusts itself to the right size by equivalence of graphs, via rescaling of the clusters that corresponds to a renaming of the labels). 

In order to properly define this stationary measure, one needs to consider a space of weighted graphs, where a given site is allowed to have a (countable) infinite family of neighbors (which corresponds to the fact that in the scaling limit, there will be infinitely many ``smaller'' clusters that are adjacent to a given macroscopic cluster).

In the case where $q >1$, things are conjecturally very similar. The only difference is that the stationary measure should be obtained by first sampling the scaling limit of the critical model, then removing a Poissonian point of pivotal points (in two dimensions, these are again ``four-arms'' points) and then defining the weighted graph in the very same way. Note that FKG comparison-type arguments show that at any large time (in the discrete dynamics, when one starts with critical weights), the configuration $w_t$ is subcritical (the cluster-size of the origin will decay exponentially). 
This leads to the following conjecture:

\medbreak
\noindent
{\bf Formulation of the universality conjecture, part one.}
{\sl For each $d \ge 2$ and all $q \ge 1$ but not too large (so that a first-order phase transition occurs), there exists a unique probability measure $\pi_{d,q}$ on  $d$-dimensional  weighted graphs (modulo equivalence) as defined above, that is invariant under the Markov process (for this value of $q$) described above: If $W_0$ is distributed according to $\pi_{d,q}$, then (up to equivalence), $W_t$ is distributed according to the same distribution at all given positive times.}

\medbreak

The second part of the universality conjecture deals with the convergence of the Markovian dynamics towards this stationary measure, when one starts with a $d$-dimensional lattice. 
Here, one needs some care in defining what it means for two (random) weighted graphs to be close. When one considers two finite weighted graphs, it is quite easy to define a distance between the two (where basically, one want the weights of the two graphs to be close in order for the two graphs to be close). Another option is to to say that two weighted finite graphs will be close if the FK model (for the appropriate value of $q$) on these two graphs are close (in total variation distance, say). Then, when we turn our attention to random $d$-dimensional graphs, there are several options to make sense of a limit (again, one wants to define this notion in such a way that one can couple nicely the FK models on the graphs); for instance, one option is to use the cut-offs $G^\eps$ and local convergence of any rooted $G^\eps$ graph. Once one has chosen such a definition, one can state the second part of the universality conjecture:

\medbreak
\noindent
{\bf Formulation of the universality conjecture, part two.}
{\sl For each $d \ge 2$ and all $q \ge 1$ but not too large (so that a first-order phase transition occurs), if $W_0=(V,c)$ is $d$-dimensional lattice with critical (for the FK$_q$ model) weights $c$, then the Markov process $W_t$ started from this lattice will converge in distribution to $\pi_{d,q}$ as $t \to \infty$. 
}

\medbreak
{Let us make another comment about the embedding:}
When one looks at the evolution of the Markov process in the stationary regime, one can define a process $(W_t, t \in \R)$ also for negative times. Each site of $V_0$ will therefore correspond to the merging of plenty of sites of $V_{-t}$ where $t$ is very large. But the $d$-dimensionality of the vertex set $V_{-t}$ suggests that one can embed in into $\R^d$ in such a way that loosely speaking, it corresponds to 
a graph that looks (on large scale) somewhat like a rather regular $d$-dimensional graph. Hence, this indicates that in the limit when $t \to \infty$, looking back in time provides a random embedding of $V_0$ in the $d$-dimensional space, that should conjecturally correspond to the actual scaling limits of clusters in the near-critical FK model on a $d$-dimensional lattice. 
In other words, from the stationary measure on these abstract equivalence classes of graphs, one can actually in fact also recover random geometric objects. Similar ideas can then be used to recover the value of critical exponents.

\medbreak
\noindent
{\bf Some results.}
There are some two-dimensional models for which the previous conjecture can be proved partially -- this corresponds to the special cases, where the scaling limit is actually understood ``geometrically'' via discrete 
holomorphicity arguments. To turn this into actual results within our setting is however far from trivial. 

\begin {itemize}
\item
As we have already pointed out, the scaling limit of critical percolation on the triangular lattice, and of the near-critical behavior is now well-understood \cite {GPS1,GPS2,GPS3}. It follows rather swiftly 
from these papers that in the case of $d=2$ and $q=1$, a stationary measure does indeed exist for the Markov process (associated to $q=1$). Uniqueness of $\pi_{2,1}$ is still an open problem.
\item 
The uniform spanning tree (or uniform spanning forest, in dimension greater than $4$) -- UST/USF in short -- can be viewed as the limit of critical FK-models when $q \to 0$. Note that in this range (for $q < 1$), the above 
arguments based on the FKG inequality that we used to describe the Markov process for infinite graphs does not hold anymore. However, due the negative correlation properties (as opposed to positive correlation 
FKG type properties for $q \ge 1$) of the UST/USF, it is nevertheless possible to 
construct the Markov process in infinite volume in this UST/USF setup, see \cite {BDW}. It is known that for a large class of planar lattices, the UST scaling limit is described by SLE$_2$ and SLE$_8$. 
This suggest that it may be possible to show that the universality conjecture basically holds in this $q=0$, $d=2$ case. 
In fact, substancial knowledge about loop-erased random walks (that are the branches in the UST) and SLE$_2$ is required. In 
 our paper with St\'ephane Benoist and Laure Dumaz \cite {BDW}, we show that the universality conjectures hold, conditional on fine upcoming results by Lawler and Viklund \cite {LV} relating the SLE$_2$
curve and its natural parametrisation to loop-erased random walks and their lenghts. 

The rooting/labeling procedure is a  little  simpler in this UST/USF case. Indeed, it is possible to follow the cluster containing a given site all the way until the scaling limit (this cluster 
ie. this tree being typically of the size of the scaling window at whatever time in the dynamics (see \cite {BDW} for details).

\end {itemize}

\medbreak
\noindent
{\bf Conclusion and outlook.} 
We have presented a formalism in which it is possible to define in elementary terms the renormalization flow as a simple Markov process on the space of discrete weighted graphs, and to view a 
``near-critical'' conjectural model as the stationary distribution for this Markov process. 
One of the features of this description is that it applies in any dimension, and separates (some of) the geometrical embedding issues from the actual rather combinatorial description of the process. 
While it does not seem to really provide obvious new avenues to prove these big universality conjectures, some related questions may be worthwhile investigating:  

\begin {itemize} 
 \item It seems to be possible to adapt the set-up to non-integer dimensions, which sheds some new light on the considerations of physicists who succesfully used dimensional expansions (for instance for $d=4- \epsilon$, at $\epsilon=1$ to obtain numerical values of exponents in dimension $3$). 
 \item Some high-dimensional results (i.e. discontinuity of the phase-transition) and maybe quantum gravity versions (\`a la Duplantier-Miller-Sheffield, possibly interpreting the weights as quantum lengths) may be accessible.
 \item Is there a way to study random planar maps with random weights in a way that generalizes the combinatorial enumerative approach and sheds some light on the stationary distributions in dimension $2$ (this may be related to the previous item)?
\item One could try to provide compactness-type arguments in order to prove abstractly the existence of some interesting measures related to these Markovian dynamics, for instance for UST models in three dimensions.
\item  It is now known that for a large class of two-dimensional lattices, the scaling limit of critical FK iterfaces when $q=2$ are described via SLE$_{16/3}$ curves. This gives hope that it may be actually
possible to 
derive (at least partially) the universality conjecture for $d=2$ and $q=2$. 
\end {itemize}

\medbreak

\noindent
{\bf Acknowledgements.}
Support and/or hospitality of the following grants and institutions is acknowledged: The Einstein Foundation Berlin, SNF-155922, NCCR Swissmap and the Isaac Newton Institute.
The present paper would probably not exist without the stimulating discussions with St\'ephane Benoist and Laure Dumaz on the joint related project \cite {BDW}.


\begin{thebibliography}{99}

   \bibitem {BDW} 
   {St\'ephane Benoist, Laure Dumaz and Wendelin Werner. 
   A renormalization approach to spanning trees.
   Preprint, 2015.} 
   
   \bibitem {Ca}
   John Cardy.
   Scaling and renormalization in statistical physics.
   CUP, 1996.
   
   \bibitem {Cetal}
Dmitry Chelkak, Hugo Duminil-Copin, Cl\'ement Hongler, Antti Kemppainen and Stanislav Smirnov.
Convergence of Ising interfaces to Schramm's SLE curves,
C. R. Math. Acad. Sci. 352 (2), 157-161 (2014)
   
\bibitem {DGP}
    Hugo Duminil-Copin, Christophe Garban and  G\'abor Pete.
The Near-Critical Planar FK-Ising Model, Communications in Mathematical Physics, 326, 1-35 (2014)

\bibitem {FK}
Cees M. Fortuin and Piet W. Kasteleyn. On the random cluster model
 I. Introduction
and relation to other models,
Physica 57, 536-564 (1972)


\bibitem{GPS1}
Christophe Garban, G\'{a}bor Pete, and Oded Schramm.
\newblock Pivotal, cluster and interface measures for critical planar percolation,
J. Amer. Math. Soc. 26, 939--1024 (2013)

\bibitem{GPS2}
Christophe Garban, G\'{a}bor Pete, and Oded Schramm.
\newblock The scaling limits of near-critical and dynamical percolation,
preprint, 2013.



\bibitem{GPS3}
Christophe Garban, G\'{a}bor Pete, and Oded Schramm.
\newblock 
The scaling limits of the Minimal Spanning Tree and Invasion Percolation in the plane,
preprint, 2013.

\bibitem{GG}
{Geoffrey Grimmett. 
The random cluster model,
Springer, 2006.
}

\bibitem {Ke}
Harry Kesten, 
Scaling relations for 2-D percolation, Comm. Math. Phys. 109, 109-156  (1987).

\bibitem{LSW} 
Gregory F. Lawler, Oded Schramm, and Wendelin Werner.
Values of Brownian intersection  exponents II: Plane exponents, Acta Mathematica 187, 275-308 (2001).

\bibitem {LSW5}
Gregory F. Lawler, Oded Schramm, and Wendelin Werner.
One-arm exponent for critical 2D percolation. Electron. J. Probab. 7, paper no. 2 (2001).

\bibitem{LSW_LERW}
Gregory F. Lawler, Oded Schramm, and Wendelin Werner.
\newblock Conformal invariance of planar loop-erased random walks and uniform
  spanning trees.
\newblock {Ann. Probab.} 32, 939--995 (2004).

\bibitem {LV}
Gregory F. Lawler, Fredrik Viklund, in preparation.

\bibitem {N}
Pierre Nolin.
Near-critical percolation in two dimensions, Electronic Journal of Probability 13, 1562-1623 (2008).

\bibitem {NW}
Pierre Nolin and Wendelin Werner. 
Asymmetry of near-critical percolation interfaces, 
J. Amer. Math. Soc. 22, 797-819 (2008).


\bibitem{S}
Oded Schramm.
\newblock Scaling limits of loop-erased random walks and uniform spanning
  trees,
\newblock { Israel J. Math.} 118, 221--288 (2000).

\bibitem {SchSm}
 Oded Schramm and Stanislav Smirnov.
On the scaling limits of planar percolation,
 Ann. Probab. 39, 1768-1814 (2011).

    
    \bibitem {Sm}{Stanislav Smirnov.}
    Discrete complex analysis and probability, Proceedings of the International Congress of Mathematicians (ICM), Hyderabad, India, 595-621 (2010). 
  
  \bibitem {SW}
  {Stanislav Smirnov and Wendelin Werner.}
  Critical exponents for two-dimensional percolation, Math. Res. Lett. 8, 729-744 (2001).
  
  
  \bibitem {T}
  Boris Tsirelson.
  Scaling limit, noise, stability, in:  Ecole d'\'et\'e  de Probabilit\'es de Saint-Flour XXXII - 2002 Lecture Notes in Mathematics 1840,  1-106 (2004).

\bibitem {Z}
    Jean Zinn-Justin, 
    Phase transitions and renormalization group, 
    OUP, 
    2013.
    
\end{thebibliography}
\end{document}